\documentclass[preprint,12pt]{elsarticle}
\usepackage{mathrsfs}
\usepackage{amsfonts}
\usepackage{amssymb,bm}
\usepackage{amsmath}
\usepackage{amsthm}
\usepackage{bbding}
\usepackage{tikz}
\usetikzlibrary{matrix,shapes,snakes}
\allowdisplaybreaks
\newtheorem{thm}{Theorem}[section]

\newtheorem{cor}[thm]{Corollary}
\newtheorem{prop}[thm]{Proposition}
\newtheorem{rmk}[thm]{Remark}

\newproof{pf}{Proof}
\newcommand{\qedd}{\hspace*{\fill}$\Box$\medskip}   

\def\ker{\mbox{\rm{ker}}}
\def\tr{{\rm{tr}}}
\def\im{{\rm{Im}}}

\renewcommand{\vec}[1]{\bm{#1}}
\def\md{{~~(\rm{mod}\;}}

\journal{Discrete Appl. Math.}
\begin{document}

\begin{frontmatter}

\title{On constructing complete permutation polynomials over finite fields of even characteristic}

\author{Baofeng Wu\corref{cor1}}
\ead{wubaofeng@iie.ac.cn}
\author{Dongdai Lin}

\cortext[cor1]{Corresponding author. Fax: +86 13426076355.}

\address{State Key Laboratory of Information Security, Institute of Information Engineering, Chinese Academy of Sciences, Beijing 100093, China}

\begin{abstract}
In this paper, a construction of complete permutation polynomials
over finite fields of even characteristic proposed by  Tu et al.
 recently is generalized in a recursive manner. Besides, several
classes of complete permutation polynomials are derived by computing
compositional inverses of known ones.
\end{abstract}


\end{frontmatter}


\section{Introduction}\label{secintro}

Let $\mathbb{F}_{q}$ be a finite field with $q$ elements where $q$
is a prime or a  prime power. A polynomial $f(x)\in\mathbb{F}_q[x]$
is called a permutation polynomial over $\mathbb{F}_{q}$ if it can
induce a bijective map from $\mathbb{F}_{q}$ to itself, and the
polynomial $f^{-1}(x)\in\mathbb{F}_q[x]$ satisfying
\[f(f^{-1}(x))\equiv f^{-1}(f(x))\equiv x \md x^q-x),\]is called
the compositional inverse of $f(x)$. Permutation polynomials have
important applications in combinatorics, coding and cryptography,
thus constructions of them have been extensively studied (see e.g.
\cite{akbary,charpin,xdhou,zbzha,zieve}). On the other hand, for
known classes of permutation polynomials, explicitly determining
their compositional inverses also attracts a lot of attention.
However, it is generally quite difficult to obtain explicit
representations of known permutation polynomials. See
\cite{wu,wu2,wang} for some recent progresses on this topic.

A permutation polynomial $f(x)\in\mathbb{F}_q[x]$ is known as a
complete permutation polynomial (CPP) over $\mathbb{F}_{q}$ if
$f(x)+x$ can permute $\mathbb{F}_{q}$ as well. Such polynomials were
initially studied by Niederreiter and Robinson in \cite{niede2}
motivated by their work on complete mappings of groups
\cite{niede1}. In fact, complete permutation polynomials can be
related to such important combinatorial objects as orthogonal latin
squares. However, to construct large classes of them is a big
challenge, and there are rare classes of complete permutation
polynomials known. We refer to \cite{mullen,lc,yuan,akbary,tu}, for
example, for some results on this topic.

Generally speaking, it seems easier  to construct complete
permutation polynomials over finite fields of even characteristic,
since it is implied by a result of Cohen that complete permutation
polynomials over $\mathbb{F}_p$ of degree $\geq2$ does not exist for
a sufficiently large prime $p$ \cite{cohen}. Very recently, several
new classes of complete permutation polynomials over finite fields
of even characteristic were constructed by Tu et al. in \cite{tu}.
More precisely, they proposed three classes of complete permutation
monomials and a class of complete permutation trinomials. Denote by
 $\tr_s^r(\cdot)$ the relative trace function from $\mathbb{F}_{2^r}$ to
 $\mathbb{F}_{2^s}$ for any positive integers $r$ and $s$ with $s\mid
 r$. Their results can be summarized in the following two theorems.

\begin{thm}[See {\cite[Theorem~1,~Theorem~2,~Theorem~3]{tu}}]\label{cppmonimial}
For two positive integers $m$, $n$, and an element
$v\in\mathbb{F}_{2^n}^*$, the monomial $v^{-1}x^d$ is a complete
permutation polynomial over $\mathbb{F}_{2^n}$ in either of the
following three cases:\\
(1) $m\geq2$, $n=3m$, $(3,m)=1$, $\tr_{m}^{n}(v)=0$, and
$d=2^{2m}+2^m+2$;\\
(2) $m\geq3$ is odd, $n=2m$,  $\tr_{m}^{n}(\beta v)=0$ or
$\tr_{m}^{n}(\beta^2 v)=0$ where $\beta$ is a primitive 3rd root of
unity in $\mathbb{F}_{2^n}^*$, and
$d=2^{m+1}+3$;\\
(3) $m\geq3$ is odd, $n=2m$,  $v$ is a non-cubic with $v^{2^m+1}=1$,
and $d=2^{m-2}(2^m+3)$.
\end{thm}

\begin{thm}[See {\cite[Theorem~4]{tu}}]\label{cpptrinomial}
For a positive integer $m$ and an element
$v\in\mathbb{F}_{2^m}\backslash\{0,1\}$, the trinomial
\[F(x)=x^{2^{2m}+1}+x^{2^m+1}+vx\]is a complete permutation
polynomial over $\mathbb{F}_{2^{3m}}$.
\end{thm}

 In this paper, we mainly focus on the class of complete permutation
 polynomials in Theorem \ref{cpptrinomial}. After noticing that the polynomial $F(x)$ in Theorem
 \ref{cpptrinomial} is just $F(x)=x\left(\tr_m^{3m}(x)+x\right)+vx$,
 we find it can be easily derived  that the polynomial
 \[\bar{F}(x)=x\left(\tr_m^{nm}(x)+x\right)+vx\]
is  a complete permutation polynomial over $\mathbb{F}_{2^{nm}}$ for
any $v\in\mathbb{F}_{2^m}\backslash\{0,1\}$ if $n$ is an odd
positive integer. Motivated by this fact, we generally consider
polynomials of the form $xL(x)+vx$, where $L(x)$ is a linearized
polynomial \cite{lidl}. Our main observation  is that complete
permutation polynomials of this form over  finite fields of
characteristic 2 can be constructed recursively. More precisely, we
find  a complete permutation polynomial of this form  over a finite
field of characteristic 2 can be obtained from a complete
permutation polynomial of the same form over certain subfield by
virtue of the relative trace function.

On the other hand, it can be easily proved that the compositional
inverse of a complete permutation polynomial also plays as a
complete one; thus new classes of complete permutation polynomials
can be derived from known ones, say, the classes of complete
permutation monomials presented in Theorem \ref{cppmonimial}, by
computing their compositional inverses. For the class of complete
permutation polynomials constructed recursively in this paper, we
can also explicitly determine the compositional inverse class thanks
to a technique given by the first author and Liu in \cite{wu},
obtaining another recursive class of complete permutation
polynomials over finite fields of even characteristic.

The rest of the paper is organized as follows. In Section
\ref{secrecucpp}, we construct a class of complete permutation
polynomials recursively to generalize Theorem \ref{cpptrinomial}. In
Section \ref{secinvcpp}, we derive several new classes of complete
permutation polynomials by computing compositional inverses of known
ones. Concluding remarks are given in Section \ref{seccon}.


\section{A construction  of CPP's generalizing Theorem
\ref{cpptrinomial}}\label{secrecucpp}

Let $m$ and $n$ be two positive integers and $q=2^m$. For
simplicity, we denote by ``$\tr$" the trace function from
$\mathbb{F}_{q^n}$ to $\mathbb{F}_q$ in the remainder of the paper.
Now we give a construction of complete permutation polynomials over
$\mathbb{F}_{q^n}$ based on complete permutation polynomials over
$\mathbb{F}_{q}$.

\begin{thm}\label{recuconst}
Let $m$ and $n$ be two positive integers where $n$ is odd, and
$q=2^m$. Assume $L(x)$ is a linearized polynomial over
$\mathbb{F}_q$ (i.e., $L(x)$ is of the form
$\sum_{i=0}^{m-1}a_ix^{2^i}$ with $a_i\in\mathbb{F}_q$, $0\leq i\leq
m-1$) such that $xL(x)+vx$ is a complete permutation polynomial over
$\mathbb{F}_q$ for some $v\in\mathbb{F}_q\backslash\{0,1\}$. Then
\[F(x)=x\left(L(\tr(x))+u\tr(x)+ux\right)+vx\]
is a complete permutation polynomial over $\mathbb{F}_{q^n}$ for any
$u\in\mathbb{F}_q$.
\end{thm}
\begin{pf}
We need only to prove that $F(x)$ can permute $\mathbb{F}_{q^n}$ for
any $u\in\mathbb{F}_q$ if $xL(x)+vx$ can permute $\mathbb{F}_q$, for
some $v\in\mathbb{F}_q\backslash\{0,1\}$. Assume $F(x)=F(y)$ for two
distinct elements $x$ and $y$ in $\mathbb{F}_{q^n}$. Since
\begin{eqnarray*}
   \tr(F(x))&=&\tr(x)L(\tr(x))+u\tr(x)^2+u\tr(x^2)+v\tr(x)  \\
   &=&\tr(x)L(\tr(x))+v\tr(x)
\end{eqnarray*}
due to the relation $\tr(x)^2=\tr(x^2)$, we have
\[\tr(x)L(\tr(x))+v\tr(x)=\tr(y)L(\tr(y))+v\tr(y),\]
which implies $\tr(x)=\tr(y)$ because $xL(x)+vx$ is a permutation
polynomial of $\mathbb{F}_q$. Then from $F(x)=F(y)$ we can get
\[(x+y)\left(L(\tr(x))+u\tr(x)\right)+u(x+y)^2=v(x+y),\]
and thus
\[L(\tr(x))+u\tr(x)+v=u(x+y)\]
as $x\neq y$. Therefore, we have
\[L(\tr(x))+u\tr(x)+v=u\tr(x+y)=0,\]
which implies $L(\tr(x))+u\tr(x)=v$. When $u\neq 0$, this leads to a
contradiction since $x^2=y^2$ follows from $F(x)=F(y)$. When $u=0$,
we have $L(\tr(x))=v\neq0$, which implies $\tr(x)\neq0$. Namely,
there exists $z\in\mathbb{F}_q^*$ such that $L(z)=v$. However, this
contradicts the fact that $x(L(x)+v)$ permutes $\mathbb{F}_q$. The
proof is completed.\qedd
\end{pf}

\begin{rmk}
It is natural to assume $v\in\mathbb{F}_q\backslash\{0,1\}$ in
Theorem \ref{recuconst} since it is necessary for $xL(x)+vx$ to be a
complete permutation polynomial over $\mathbb{F}_q$. In fact,
$xL(x)$ and $xL(x)+x$ cannot be permutation polynomials
simultaneously. Indeed, $L(x)$ is necessarily a linearized
permutation polynomial, which implies $1\in\im(L)$, if $xL(x)$ is a
permutation polynomial, while $1\not\in\im(L)$ if $xL(x)+x$ is a
permutation polynomial, where $\im(L)$ represents the image space of
the linear transformation induced by $L(x)$.
\end{rmk}

It is easy to see that the complete permutation polynomial $F(x)$
constructed in Theorem \ref{recuconst} is still of the form
$x\bar{L}(x)+vx$ for certain linearized polynomial $\bar{L}(x)$ over
$\mathbb{F}_{q^n}$. Therefore, Theorem \ref{recuconst} actually says
that complete permutation polynomials of the form  $x\bar{L}(x)+vx$
over a finite field of characteristic 2 can be used to construct
complete permutation polynomials of the same form over odd-degree
extensions of this field; thus it presents a recursive construction
of complete permutation polynomials.

By fixing $L(x)=0$ in Theorem \ref{recuconst}, we can obtain the
following construction of complete permutation polynomials which
generalizes the one given in Theorem \ref{cpptrinomial}.

\begin{cor}\label{cpptrigeneralize}
For a positive integer $m$, an odd positive integer $n$ and any
$u\in\mathbb{F}_{2^m}$, $v\in\mathbb{F}_{2^m}\backslash\{0,1\}$, the
polynomial
\[F(x)=x\left(u\tr_m^{nm}(x)+ux\right)+vx\]is a complete permutation
polynomial over $\mathbb{F}_{2^{nm}}$.
\end{cor}

Combining Theorem \ref{recuconst} and Corollary
\ref{cpptrigeneralize}, we can also obtain the following
construction of complete permutation polynomials involving
multi-trace terms.

\begin{cor}\label{cpptwistgen}
Let $m$ and $n$ be two positive integers where $n$ is odd, and
$q=2^m$. Assume $d_1,~d_2,\ldots,d_s$ are distinct positive integers
with $d_1\mid d_2\mid \cdots\mid d_s\mid n$ and let
$c_0\in\mathbb{F}_q$, $c_i\in\mathbb{F}_{q^{d_i}}$, $1\leq i\leq s$.
Then for any $\tilde{c}\in\mathbb{F}_q\backslash\{0,1\}$, the
polynomial
\[F(x)=x\left(\sum_{i=0}^{c_i}\tr_{d_im}^{nm}(x)+cx\right)+\tilde{c}x\]
 is a complete permutation polynomial over
$\mathbb{F}_{q^n}$, where $c=\sum_{j=0}^sc_j$ and $d_0=1$.
\end{cor}
\begin{pf}
Let $L_1(x)=c_0\tr_{d_0m}^{d_1m}(x)+c_0x$ and for $2\leq i\leq s+1$,
let
\[L_i(x)=L_{i-1}\left(\tr^{d_im}_{d_{i-1}m}(x)\right)
+\left(\sum_{j=0}^{i-1}c_j\right)\tr^{d_im}_{d_{i-1}m}(x)+\left(\sum_{j=0}^{i-1}c_j\right)x,\]
where $d_{s+1}=n$. Since $F_1(x)=xL_1(x)+\tilde{c}x$ is a complete
permutation polynomial over $\mathbb{F}_{q^{d_1}}$ according to
Corollary \ref{cpptrigeneralize}, we know by induction that
\[F_i(x)=xL_i(x)+\tilde{c}x\]
is a complete permutation polynomial over $\mathbb{F}_{q^{d_i}}$ for
any $2\leq i\leq s+1$  from Theorem \ref{recuconst}. According to
the transitivity of the trace function, it can be easily derived by
induction  that
\[L_i(x)=\sum_{j=0}^{i-1}c_j\tr^{d_im}_{d_{j}m}(x)+\left(\sum_{j=0}^{i-1}c_j\right)x\]
 for $2\leq i\leq s+1$. Hence finally we know
that $F(x)=F_{s+1}(x)$ is a complete permutation polynomial over
$\mathbb{F}_{q^n}$.\qedd
\end{pf}


\section{Constructing CPP's by inverting known
ones}\label{secinvcpp}

In this section, we propose another approach to obtain complete
permutation polynomials. The main observation is included in the
following proposition.

\begin{prop}\label{invcpp}
Let $f(x)$ be a complete permutation polynomial over $\mathbb{F}_Q$
where $Q$ is a primer power. Then $f^{-1}(x)$ is also a complete
permutation polynomial over $\mathbb{F}_{Q}$.
\end{prop}
\begin{pf}
Assume $g(x)=f^{-1}(x)+x$. Then we have
\[g(f(x))=f^{-1}(f(x))+f(x)=x+f(x).\]
Since $x+f(x)$ and $f(x)$ are both permutation polynomials over
$\mathbb{F}_{Q}$, we know that $g(x)$ is also a permutation
polynomial over $\mathbb{F}_{Q}$, which implies  $f^{-1}(x)$ is a
complete permutation polynomial over $\mathbb{F}_{Q}$. \qedd
\end{pf}

By Proposition \ref{invcpp}, we can get complete permutation
polynomials via computing compositional inverses of known ones. But
generally speaking, it is  far from a simple matter to obtain
explicit representations of compositional inverses of known
permutation polynomials over finite fields. However, in some special
cases, say, the permutation polynomials in consideration are
permutation monomials, there are less difficulties, of course, to
compute compositional inverses. This is because $ax^d$ is a
permutation polynomial over a finite field $\mathbb{F}_Q$ if and
only if $(d,Q-1)=1$ and when this condition holds, the compositional
inverse of $ax^d$ is just $(x/a)^{d^{-1}}$, where $d^{-1}$ is a
positive integer satisfying $dd^{-1}\equiv1\md Q-1)$. Thus to
explicitly get compositional inverses of permutation monomials over
$\mathbb{F}_Q$ is equivalent to explicitly get inverses of elements
in $(\mathbb{Z}/(Q-1)\mathbb{Z})^*$, the unit group of the integer
residue ring $\mathbb{Z}/(Q-1)\mathbb{Z}$, which can be done by the
Euclid algorithm.

Consequently, new classes of complete permutation monomials in the
following three theorems can be directly  derived from those given
in Theorem~1, Theorem~2 and Theorem~3 of \cite{tu}, respectively, by
computing their compositional inverses. Proofs of Theorem
\ref{monomialthm1}, Theorem \ref{monomialeg} and Theorem
\ref{monomialthm2} can be directly obtained, which will be omitted
here. Note that for any $a\in\mathbb{F}_Q^*$, $af(x/a)$ is a
complete permutation polynomial over  $\mathbb{F}_Q$ if and only if
$f(x)$ is. This fact needs to be applied in the proofs.

\begin{thm}\label{monomialthm1}
Let $m\geq 2$ be an integer with $(3,m)=1$ and
$v\in\mathbb{F}_{2^{3m}}^*$ with $\tr_m^{3m}(v)=0$. Then the
monomial $vx^{2^{3m-1}+2^{3m-2}-2^{2m-2}-2^{m-2}}$ is a complete
permutation polynomial over $\mathbb{F}_{2^{3m}}$.
\end{thm}

\begin{thm}\label{monomialeg}
Let $m\geq 3$ be an odd integer and $v\in\mathbb{F}_{2^{2m}}^*$ with
$\tr_m^{2m}(\beta v)=0$ or $\tr_m^{2m}(\beta^2 v)=0$, where $\beta$
 is a primitive 3rd root of unity in $\mathbb{F}_{2^{2m}}$.
 Then the monomial $vx^{{({2^{2m+1}-2^{m+1}+1})/{5}}}$ or
 $vx^{{({2^{2m}-2^{m+1}+2})/{5}}}$, respectively, is a complete
permutation polynomial over $\mathbb{F}_{2^{2m}}$ according as
$m\equiv1\md4)$ or $m\equiv3\md4)$, respectively.
\end{thm}

\begin{thm}\label{monomialthm2}
Let $m\geq 3$ be an odd integer and $v\in\mathbb{F}_{2^{m}}^*$ be a
non-cubic with $v^{2^m+1}=1$. Then the monomial
$vx^{2^{2m-1}+2^m+2^{m-1}-1}$ is a complete permutation polynomial
over $\mathbb{F}_{2^{2m}}$.
\end{thm}

In fact, for the exponent $d$ appearing in the monomials in  Theorem
2 or Theorem 3 of \cite{tu}, it is not so easy to explicitly get
$d^{-1}$ directly using the Euclid algorithm. To overcome the
difficulties, we can apply the Chinese Remainder Theorem. As an
example, we show how to compute $(2^{m+1}+3)^{-1}$ in
$(\mathbb{Z}/(2^{2m}-1)\mathbb{Z})^*$ for an odd integer $m\geq3$ to
obtain Theorem \ref{monomialeg}.

 For any
$r\in(\mathbb{Z}/(2^{2m}-1)\mathbb{Z})^*$, denote by $r_1^{-1}$ and
$r_2^{-1}$  the inverses of $r_1$  in
$(\mathbb{Z}/(2^{m}-1)\mathbb{Z})^*$ and $r_2$ in
$(\mathbb{Z}/(2^{m}+1)\mathbb{Z})^*$, respectively, where $r\equiv
r_1\md 2^m-1)$ and $r\equiv r_2\md 2^m+1)$. Then by the Chinese
Remainder Theorem, we know that
\begin{equation}\label{eqnCRT}
r^{-1}\equiv 2^{m-1}(2^m+1)r_1^{-1}+2^{m-1}(2^m-1)r_2^{-1}\md
2^{2m}-1).
\end{equation}
Now for $r=2^{m+1}+3$, it is obvious that $r_1=5$ and $r_2=1$; thus
$r_2^{-1}=1$. To get $r_1^{-1}$, we note that
$5^{-1}=3\cdot(2^4-1)^{-1}$ since $5={(2^4-1)}/{3}$. By the Euclid
algorithm, it is easy to derive that
\[(2^4-1)^{-1}=\left\{\begin{array}{ll}
-2\cdot\frac{2^{m-1}-1}{2^4-1}&\text{if}~m\equiv1\md4)\\
1+2^4\cdot\frac{2^{m-1}-1}{2^4-3}&\text{if}~m\equiv3\md4),
\end{array}\right.\]
which implies
\[5^{-1}=\left\{\begin{array}{ll}
-2\cdot\frac{2^{m-1}-1}{5}&\text{if}~m\equiv1\md4)\\
3+2^4\cdot\frac{2^{m-1}-1}{5}&\text{if}~m\equiv3\md4).
\end{array}\right.\]
Then we get  that
\[(2^{m+1}+3)^{-1}=\left\{\begin{array}{ll}
\frac{2^{2m+1}-2^{m+1}+1}{5}&\text{if}~m\equiv1\md4)\\
\frac{2^{2m}-2^{m+1}+2}{5}&\text{if}~m\equiv3\md4)
\end{array}\right.\]
from \eqref{eqnCRT}. Theorem \ref{monomialeg} follows.\bigskip

The main purpose of this section is to obtain complete permutation
polynomials by explicitly representing the compositional inverse of
the polynomial $F(x)$ defined in Theorem \ref{recuconst}. In the
sequel, we give the formula of the compositional inverse of $F(x)$
and verify it at first, and explain the process to derive it
afterwards. Firstly we remark that in representing polynomials over
a finite field $\mathbb{F}_{2^e}$, we sometimes use $x^{1/2}$ and
$1/x$ instead of $x^{2^{e-1}}$ and $x^{2^e-2}$, respectively.

\begin{thm}\label{recuinvconst}
Let $m$ and $n$ be two positive integers where $n$ is odd, and
$q=2^m$. Assume $L(x)$ is a linearized polynomial over
$\mathbb{F}_q$ such that $xL(x)+vx$ is a complete permutation
polynomial over $\mathbb{F}_q$ for some
$v\in\mathbb{F}_q\backslash\{0,1\}$. Let $g(x)$ be the compositional
inverse of $xL(x)+vx$ over $\mathbb{F}_q$. Then,\\
(1) the polynomial
\[\bar{F}(x)=\frac{x}{v}+\left(\frac{g(\tr(x))}{\tr(x)}+\frac{1}{v}\right)x\tr(x)^{q-1}\]
is a complete polynomial over $\mathbb{F}_{q^n}$; and \\
(2) the polynomial
\begin{equation*}
\begin{aligned}
&\tilde{F}(x)\\
&=\left(1+\tr(x)^{q-1}\right)\sum_{j=0}^{m-1}{\frac{u^{2^j-1}}{v^{2^{j+1}-1}}}
 \left(\sum_{k=0}^{\frac{n-1}{2}}x^{q^{2k}}\right)^{2^j}\\
 &\quad+\left[\tr(x)^{q-1}+\left(u^{1/2}g(\tr(x))\tr(x)+\tr(x)^{3/2}\right)^{q-1}\right]\left(\frac{x}{u}\right)^{1/2}\\
&\quad+\left(u^{1/2}g(\tr(x))\tr(x)+\tr(x)^{3/2}\right)^{q-1}\\
&\times\left[g\left(\tr(x)\right)+u^{1/2}\sum_{j=0}^{m-1}{\left(\frac{\tr(x)}{u^{1/2}g\left(\tr(x)\right)}+u^{1/2}g\left(\tr(x)\right)\right)^{-(2^{j+1}-1)}}
\left(\sum_{k=0}^{\frac{n-1}{2}}x^{q^{2k}}\right)^{2^j}\right]
\end{aligned}
\end{equation*}
is a complete permutation polynomial over $\mathbb{F}_{q^n}$ for any
$u\in\mathbb{F}_q^*$.
\end{thm}
\begin{pf}
We proceed by directly verify that $\bar{F}(F(x))=x$ when $u=0$ and
$\tilde{F}(F(x))=x$ when $u\in\mathbb{F}_q^*$ for any
$x\in\mathbb{F}_{q^n}$, then the result follows from Proposition
\ref{invcpp}.

First, note that for any $u\in\mathbb{F}_q$,
\[\tr({F}(x))=\tr(x)L(\tr(x))+v\tr(x),\] which implies
\[g(\tr(F(x)))=g\left(\tr(x)L(\tr(x))+v\tr(x)\right)=\tr(x)\]since
$g\left(yL(y)+vy\right)=y$ for any $y\in\mathbb{F}_q$.

(1) When $u=0$, we have $F(x)=xL(\tr(x))+vx$. Note that
$\tr(F(x))=0$ if and only if $\tr(x)=0$.  Thus for
$x\in\mathbb{F}_{q^n}$ with $\tr(x)=0$, we have
\[\bar{F}(F(x))=\frac{F(x)}{v}=\frac{vx}{v}=x,\]and for $x\in\mathbb{F}_{q^n}$ with
$\tr(x)\neq0$, we have
\[\bar{F}(F(x))=\frac{F(x)g(\tr(F(x)))}{\tr(F(x))}=\frac{\left(xL(\tr(x))+vx\right)\tr(x)}{\tr(x)L(\tr(x))+v\tr(x)}=x.\]
Hence $\bar{F}(F(x))=x$ for any $x\in\mathbb{F}_{q^n}$.

(2) When $u\in\mathbb{F}_q^*$, it is easy to see that
\[\tilde{F}(x)=\sum_{j=0}^{m-1}{\frac{u^{2^j-1}}{v^{2^{j+1}-1}}}
 \sum_{k=0}^{\frac{n-1}{2}}x^{2^{2km+j}}\]
for $x\in\mathbb{F}_{q^n}$ with $\tr(x)=0$,
\[\tilde{F}(x)=\left(\frac{x}{u}\right)^{1/2}\] for $x\in\mathbb{F}_{q^n}$ with $\tr(x)\neq0$ and
$u^{1/2}g(\tr(x))+\tr(x)^{1/2}=0$, and
 \[\tilde{F}(x)=g\left(\tr(x)\right)+u^{-1/2}\sum_{j=0}^{m-1}{\left(\frac{\tr(x)}{u^{1/2}g(\tr(x))}+u^{1/2}g(\tr(x))\right)^{-(2^{j+1}-1)}}
 \sum_{k=0}^{\frac{n-1}{2}}x^{2^{2km+j}}\]otherwise.

Recall that $\tr(F(x))=0$ if and only if $\tr(x)=0$.  For
$x\in\mathbb{F}_{q^n}$ with $\tr(x)=0$, we have $F(x)=ux^2+vx$; thus
\begin{eqnarray*}
  \tilde{F}(F(x)) &=&\sum_{j=0}^{m-1}{\frac{u^{2^j-1}}{v^{2^{j+1}-1}}}
 \sum_{k=0}^{\frac{n-1}{2}}(ux^2+vx)^{2^{2km+j}}  \\
   &=&\sum_{j=0}^{m-1}{\frac{u^{2^{j+1}-1}}{v^{2^{j+1}-1}}}
 \sum_{k=0}^{\frac{n-1}{2}}x^{2^{2km+j+1}}+
 \sum_{j=0}^{m-1}{\frac{u^{2^j-1}}{v^{2^{j}-1}}}
 \sum_{k=0}^{\frac{n-1}{2}}x^{2^{2km+j}}\\
 &=&\sum_{k=0}^{\frac{n-1}{2}}x^{q^{2k+1}}+\sum_{k=0}^{\frac{n-1}{2}}x^{q^{2k}}\\
 &=&\tr(x)+x\\
 &=&x.
\end{eqnarray*}
On the other hand, for $x\in\mathbb{F}_{q^n}$ with $\tr(x)\neq0$, we
have
\[u^{1/2}g\left(\tr(F(x))\right)+\tr(F(x))^{1/2}=\tr(x)^{1/2}\left(L(\tr(x))+u\tr(x)+v\right)^{1/2},\]
which implies $u^{1/2}g\left(\tr(F(x))\right)+\tr(F(x))^{1/2}=0$ if
and only if $L(\tr(x))+u\tr(x)+v=0$. Therefore, when
$L(\tr(x))+u\tr(x)+v=0$, we have $F(x)=ux^2$, and thus
\[\tilde{F}(F(x))=\left(\frac{F(x)}{u}\right)^{1/2}=\left(\frac{ux^2}{u}\right)^{1/2}=x;\]
when $L(\tr(x))+u\tr(x)+v\neq0$, we have
\begin{eqnarray*}
   \tilde{F}(F(x))&=&g\left(\tr(F(x))\right)+u^{-1/2}\sum_{j=0}^{m-1}{\left(\frac{\tr(F(x))}{u^{1/2}g(\tr(F(x)))}+u^{1/2}g(\tr(F(x)))\right)^{-(2^{j+1}-1)}}\\
 &&\qquad\qquad\qquad\qquad\times\sum_{k=0}^{\frac{n-1}{2}}F(x)^{2^{2km+j}}\\
   &=&\tr(x)+u^{-1/2}\sum_{j=0}^{m-1}\frac{u^{2^j-1/2}}{\left(L(\tr(x))+u\tr(x)+v\right)^{2^{j+1}-1}}
 \\&&\qquad\quad\qquad\times\sum_{k=0}^{\frac{n-1}{2}}\left[x\left(L(\tr(x))+u\tr(x)+ux\right)+vx\right]^{2^{2km+j}}  \\
   &=&\tr(x)+\sum_{j=0}^{m-1}\frac{\left(L(\tr(x))+u\tr(x)+v\right)^{2^{j}}u^{2^j-1}}{\left(L(\tr(x))+u\tr(x)+v\right)^{2^{j+1}-1}}
 \sum_{k=0}^{\frac{n-1}{2}}x^{2^{2km+j}}  \\
 &&\qquad~+\sum_{j=0}^{m-1}\frac{u^{2^{j+1}-1}}{\left(L(\tr(x))+u\tr(x)+v\right)^{2^{j+1}-1}}
 \sum_{k=0}^{\frac{n-1}{2}}x^{2^{2km+j+1}}\\
 &=&\tr(x)+\sum_{j=0}^{m-1}\frac{u^{2^j-1}}{\left(L(\tr(x))+u\tr(x)+v\right)^{2^{j}-1}}
 \sum_{k=0}^{\frac{n-1}{2}}x^{2^{2km+j}}\\&&\qquad~+\sum_{j=0}^{m-1}\frac{u^{2^{j+1}-1}}{\left(L(\tr(x))+u\tr(x)+v\right)^{2^{j+1}-1}}
 \sum_{k=0}^{\frac{n-1}{2}}x^{2^{2km+j+1}}\\
 &=&\tr(x)+\sum_{k=0}^{\frac{n-1}{2}}x^{q^{2k}}+\sum_{k=0}^{\frac{n-1}{2}}x^{q^{2k+1}}\\
 &=&x.
\end{eqnarray*}
Finally we know that $\tilde{F}(F(x))=x$ for any
$x\in\mathbb{F}_{q^n}$. This completes the proof. \qedd
\end{pf}

The polynomials $\bar{F}(x)$ and $\tilde{F}(x)$ in Theorem
\ref{recuinvconst}, though having complicated representations, are
indeed complete permutation polynomials over $\mathbb{F}_{q^n}$
induced from a complete permutation polynomial $g(x)$ over
$\mathbb{F}_{q}$. Hence Theorem \ref{recuinvconst} also provides
recursive constructions of complete permutation polynomials over
odd-degree extensions of finite fields.

 By fixing $L(x)=0$ in Theorem
\ref{recuinvconst}, we can obtain the following class of complete
permutation polynomials. They correspond to the class of complete
permutation polynomials proposed in Corollary
\ref{cpptrigeneralize}.

\begin{cor}
For a positive integer $m$, an odd positive integer $n$ and any
$u\in\mathbb{F}_{2^m}^*$, $v\in\mathbb{F}_{2^m}\backslash\{0,1\}$,
the polynomial
\begin{eqnarray*}
\tilde{F}(x)
&=&\left(1+\tr_{m}^{nm}(x)^{q-1}\right)\sum_{j=0}^{m-1}{\frac{u^{2^j-1}}{v^{2^{j+1}-1}}}
 \left(\sum_{k=0}^{\frac{n-1}{2}}x^{q^{2k}}\right)^{2^j}\\
 &&+\left[\tr_{m}^{nm}(x)^{q-1}+\left(u\tr_{m}^{nm}(x^2)+v^2\tr_{m}^{nm}(x)\right)^{q-1}\right]\left(\frac{x}{u}\right)^{1/2}\\
&&+\left(u\tr_{m}^{nm}(x^2)+v^2\tr_{m}^{nm}(x)\right)^{q-1}\\
&&\times\left[\frac{\tr_{m}^{nm}(x)}{v}+u^{1/2}\sum_{j=0}^{m-1}{\left(\frac{u\tr_{m}^{nm}(x)+v^2}{u^{1/2}v}\right)^{-(2^{j+1}-1)}}
\left(\sum_{k=0}^{\frac{n-1}{2}}x^{q^{2k}}\right)^{2^j}\right]
\end{eqnarray*}
is a complete permutation polynomial over $\mathbb{F}_{2^{nm}}$.
\end{cor}
\begin{pf}
Note that when $L(x)=0$ in Theorem \ref{recuinvconst}, we have
$g(x)=x/v$, and thus
\[\left(u^{1/2}g(\tr_{m}^{nm}(x))\tr_{m}^{nm}(x)+\tr_{m}^{nm}(x)^{3/2}\right)^{q-1}=\left(u\tr_{m}^{nm}(x^2)+v^2\tr_{m}^{nm}(x)\right)^{q-1}.\]
Then the result follows from Theorem \ref{recuinvconst}.\qedd
\end{pf}

In addition, according to the proof of Corollary \ref{cpptwistgen},
compositional inverses of complete permutation polynomials from
Corollary \ref{cpptwistgen} can also be explicitly obtained from
Theorem \ref{recuinvconst} by induction, which will lead to another
class of complete permutation polynomials involving multi-trace
terms. However, their representations are rather complicated and the
derivation of them will be left to the interested readers.

For completeness of this paper, in the following we explain how to
compute the compositional inverse of the permutation polynomial
$F(x)$ given in Theorem \ref{recuconst}
 for any $u\in\mathbb{F}_q$. The main idea has already been indicated in
 \cite{wu}. Consider a graph of maps
\begin{center}
\begin{tikzpicture}
\matrix (m) [matrix of math nodes,row sep=4em,column sep=5em,minimum
width=2em]
  {
     \mathbb{F}_{q^n} & \mathbb{F}_{q^n} \\
    \mathbb{F}_{q}\oplus\ker(\tr) & \mathbb{F}_{q}\oplus\ker(\tr) \\};
  \path[-stealth]
    (m-1-1) edge node [left] {$\phi$} (m-2-1)
            edge node [above] {$F(x)$} (m-1-2)
    (m-2-1) edge node [above] {$\vec F(y,z)$} (m-2-2)
    (m-1-2) edge node [right] {$\phi$} (m-2-2);
\end{tikzpicture}
\end{center}
Here the map $\phi$, which is defined by
\[\phi(x)=(\tr(x),~x+\tr(x)),\]
can induce an isomorphism between $\mathbb{F}_{q^n}$ and
$\mathbb{F}_{q}\oplus\ker(\tr)$, where $\ker(\tr)$ represents the
kernel space of the trace map of $\mathbb{F}_{q^n}$ over
$\mathbb{F}_{q}$, because $n$ is odd, and $\vec F(y,z)$ is a
bivariate polynomial system that can make the above graph
commutative (i.e., induce a bijective map from
$\mathbb{F}_{q}\oplus\ker(\tr)$ to itself). To obtain $F^{-1}(x)$
(the compositional inverse of $F(x)$), we need only to find the
polynomial system $\vec F^{-1}(y,z)$ that can induce the inverse map
of the map induced by $\vec F(y,z)$ on
$\mathbb{F}_{q}\oplus\ker(\tr)$.

For any $x\in\mathbb{F}_{q^n}$, we let $y=\tr(x)$ and $z=x+\tr(x)$.
Since
\[\tr(F(x))=\tr(x)L(\tr(x))+v\tr(x)=yL(y)+vy\]and
\begin{eqnarray*}
  F(x)+\tr(F(x)) &=&(y+z)L(y)+u(y+z)y+u(y+z)^2+v(y+z)+yL(y)+vy\\
   &=&uz^2+\left(L(y)+uy+v\right)z,
\end{eqnarray*}
we have
\[\vec F(y,z)=\left(yL(y)+vy,~z^2+\left(L(y)+uy+v\right)z\right).\]
Assume $F(x)=X$ and $Y=\tr(X)$, $Z=X+\tr(X)$. Consider the system of
equations
\begin{equation}\label{eqnsyst}
\left\{\begin{aligned}
&yL(y)+vy=Y\\
&uz^2+\left(L(y)+uy+v\right)z=Z .
\end{aligned}\right.
\end{equation}
It is direct to obtain that $y=g(Y)$ where $g$ is defined in Theorem
\ref{recuinvconst}. Thus  $L(y)=L(g(Y))=0$ when $Y=0$ and
$L(y)=L(g(Y))=Y/g(Y)+v$ when $Y\neq0$ since
\[g(Y)L(g(Y))+vg(Y)=Y.\]

First we consider the case for $u=0$. When $Y=0$, we can get from
\eqref{eqnsyst}  that $y=0$ and $z=Z/v$, i.e.,
\[\vec F^{-1}(Y,Z)=\left(0,\frac{Z}{v}\right).\]
Hence when $\tr(X)=0$ we have
\[F^{-1}(X)=\frac{Z}{v}=\frac{X+\tr(X)}{v}=\frac{X}{v}.\]
When $Y\neq0$, we can get from \eqref{eqnsyst}  that
\[\vec F^{-1}(Y,Z)=\left(g\left(Y\right),\frac{Z}{L(g(Y))+v}\right)=\left(g\left(Y\right),\frac{Zg(Y)}{Y}\right).\]
Hence when $\tr(X)\neq0$ we have
\[F^{-1}(X)=g\left(Y\right)+\frac{Zg(Y)}{Y}=g(\tr(X))+\frac{X+\tr(X)}{\tr(X)}g(\tr(X))=\frac{Xg(\tr(X))}{\tr(X)}.\]
Therefore, in the case $u=0$, by Lagrange interpolation we have
\begin{eqnarray*}
F^{-1}(X)&=&\frac{X}{v}\left(1+\tr(X)^{q-1}\right) + \frac{Xg(\tr(X))}{\tr(X)}\tr(X)^{q-1}\\
 &=&\frac{X}{v}+\left(\frac{g(\tr(X))}{\tr(X)}+\frac{1}{v}\right)X\tr(X)^{q-1}.
\end{eqnarray*}
This yields $\bar{F}(x)$ in Theorem \ref{recuinvconst}.

Now we consider the case for $u\in\mathbb{F}_q^*$. When $Y=0$, we
can get from \eqref{eqnsyst}  that $y=0$ and $z^2+vz/u=Z/u$.
According to \cite[Lemma~3.3]{wu}, we have
\[z=\sum_{j=0}^{m-1}{\left(\frac{v}{u}\right)^{-(2^{j+1}-1)}}
 \sum_{k=0}^{\frac{n-1}{2}}\left(\frac{Z}{u}\right)^{2^{2km+j}}
 =\sum_{j=0}^{m-1}{\frac{u^{2^j-1}}{v^{2^{j+1}-1}}}
 \sum_{k=0}^{\frac{n-1}{2}}Z^{2^{2km+j}}.\]
Therefore, in this case,
\[\vec F^{-1}(Y,Z)=\left(0,\sum_{j=0}^{m-1}{\frac{u^{2^j-1}}{v^{2^{j+1}-1}}}
 \sum_{k=0}^{\frac{n-1}{2}}Z^{2^{2km+j}}\right).\]
Hence when $\tr(X)=0$ we have
\begin{eqnarray}\label{F1}
 F^{-1}(X)  &=&\sum_{j=0}^{m-1}{\frac{u^{2^j-1}}{v^{2^{j+1}-1}}}
 \sum_{k=0}^{\frac{n-1}{2}}(X+\tr(X))^{2^{2km+j}}\nonumber\\
   &=&\sum_{j=0}^{m-1}{\frac{u^{2^j-1}}{v^{2^{j+1}-1}}}
 \left(\sum_{k=0}^{\frac{n-1}{2}}X^{q^{2k}}\right)^{2^j}.
\end{eqnarray}

When $Y\neq0$, we can get from \eqref{eqnsyst}  that $y=g(Y)$ and
$$z^2+\left(\frac{L(y)+v}{u}+y\right)z=z^2+\left(\frac{Y}{ug(Y)}+g(Y)\right)z=Z/u.$$
We distinguish two cases. First, when ${Y}/{(ug(Y))}+g(Y)=0$, i.e.,
$u^{1/2}g(Y)+Y^{1/2}=0$, we have $z=(Z/u)^{1/2}$. This implies
\[\vec F^{-1}(Y,Z)=\left(g\left(Y\right),\left(\frac{Z}{u}\right)^{1/2}\right)
=\left(\left(\frac{Y}{u}\right)^{1/2},\left(\frac{Z}{u}\right)^{1/2}\right).\]
Hence when $\tr(X)\neq0$ and $u^{1/2}g(\tr(X))+\tr(X)^{1/2}=0$ we
have
\begin{equation}\label{F2}
F^{-1}(X)=\left(\frac{Y}{u}\right)^{1/2}+\left(\frac{Z}{u}\right)^{1/2}=\left(\frac{X}{u}\right)^{1/2}.
\end{equation}
Second, when ${Y}/{(ug(Y))}+g(Y)\neq0$, i.e.,
$u^{1/2}g(Y)+Y^{1/2}\neq0$, according to \cite[Lemma~3.3]{wu} we
have
\begin{eqnarray*}
   z&=&\sum_{j=0}^{m-1}{\left(\frac{Y}{ug\left(Y\right)}+g\left(Y\right)\right)^{-(2^{j+1}-1)}}
 \sum_{k=0}^{\frac{n-1}{2}}\left(\frac{Z}{u}\right)^{2^{2km+j}}  \\
   &=&u^{1/2}\sum_{j=0}^{m-1}{\left(\frac{Y}{u^{1/2}g\left(Y\right)}+u^{1/2}g\left(Y\right)\right)^{-(2^{j+1}-1)}}
 \sum_{k=0}^{\frac{n-1}{2}}Z^{2^{2km+j}}.
\end{eqnarray*}
This implies
\[\vec F^{-1}(Y,Z)=\left(g\left(Y\right),u^{1/2}\sum_{j=0}^{m-1}{\left(\frac{Y}{u^{1/2}g\left(Y\right)}+u^{1/2}g\left(Y\right)\right)^{-(2^{j+1}-1)}}
 \sum_{k=0}^{\frac{n-1}{2}}Z^{2^{2km+j}}\right).\]
Hence when $u^{1/2}g(\tr(X))+\tr(X)^{1/2}\neq0$ we have
\begin{equation}\label{F3}
\begin{aligned}
&F^{-1}(X)\\
  &=g\left(\tr(X)\right)+u^{1/2}\sum_{j=0}^{m-1}{\left(\frac{\tr(X)}{u^{1/2}g\left(\tr(X)\right)}+u^{1/2}g\left(\tr(X)\right)\right)^{-(2^{j+1}-1)}}\\
  &\qquad\quad\qquad~\qquad~\times \sum_{k=0}^{\frac{n-1}{2}}(X+\tr(X))^{2^{2km+j}} \\
   &=g\left(\tr(X)\right)+u^{1/2}\sum_{j=0}^{m-1}{\left(\frac{\tr(X)}{u^{1/2}g\left(\tr(X)\right)}+u^{1/2}g\left(\tr(X)\right)\right)^{-(2^{j+1}-1)}}
\sum_{k=0}^{\frac{n-1}{2}}X^{2^{2km+j}}  \\
&\quad+\frac{n+1}{2}u^{1/2}\sum_{j=0}^{m-1}{\left(\frac{\tr(X)}{u^{1/2}g\left(\tr(X)\right)}+u^{1/2}g\left(\tr(X)\right)\right)^{-(2^{j+1}-1)}}\tr(X)^{2^j}\\
&=g\left(\tr(X)\right)+u^{1/2}\sum_{j=0}^{m-1}{\left(\frac{\tr(X)}{u^{1/2}g\left(\tr(X)\right)}+u^{1/2}g\left(\tr(X)\right)\right)^{-(2^{j+1}-1)}}
\left(\sum_{k=0}^{\frac{n-1}{2}}X^{q^{2k}}\right)^{2^j}
\end{aligned}
\end{equation}
because
\begin{equation*}
\begin{aligned}
&\sum_{j=0}^{m-1}{\left(\frac{\tr(X)}{u^{1/2}g\left(\tr(X)\right)}+u^{1/2}g\left(\tr(X)\right)\right)^{-(2^{j+1}-1)}}\tr(X)^{2^j}\\
   &=\left(\frac{\tr(X)}{u^{1/2}g\left(\tr(X)\right)}+u^{1/2}g\left(\tr(X)\right)\right)
\sum_{j=0}^{m-1}\left(\frac{\tr(X)^{2^j}}{u^{2^j}g\left(\tr(X)\right)^{2^{j+1}}}+
\frac{u^{2^j}g\left(\tr(X)\right)^{2^{j+1}}}{\tr(X)^{2^j}}\right)^{-1}\\&=0
\end{aligned}
\end{equation*}
according to \cite[Lemma~3.5]{wu}.

We denote by $F_1(X)$, $F_2(X)$ and $F_3(X)$ the formulas of
$F^{-1}(X)$ in \eqref{F1}, \eqref{F2} and \eqref{F3}, respectively.
Finally, $F^{-1}(X)$  can be represented by Lagrange interpolation,
that is
\begin{equation*}
\begin{aligned}
&F^{-1}(X)\\
&=\left(1+\tr(X)^{q-1}\right)F_1(X)+\left[1+\left(u^{1/2}g(\tr(X))+\tr(X)^{1/2}\right)^{q-1}\right]\tr(X)^{q-1}F_2(X)\\
  &\quad+\left(u^{1/2}g(\tr(X))+\tr(X)^{1/2}\right)^{q-1}\tr(X)^{q-1}F_3(X)\\
&=\left(1+\tr(X)^{q-1}\right)\sum_{j=0}^{m-1}{\frac{u^{2^j-1}}{v^{2^{j+1}-1}}}
 \left(\sum_{k=0}^{\frac{n-1}{2}}X^{q^{2k}}\right)^{2^j}\\
 &\quad+\left[\tr(X)^{q-1}+\left(u^{1/2}g(\tr(X))\tr(X)+\tr(X)^{3/2}\right)^{q-1}\right]\left(\frac{X}{u}\right)^{1/2}\\
&\quad+\left(u^{1/2}g(\tr(X))\tr(X)+\tr(X)^{3/2}\right)^{q-1}\\
&\times\left[g\left(\tr(X)\right)+u^{1/2}\sum_{j=0}^{m-1}{\left(\frac{\tr(X)}{u^{1/2}g\left(\tr(X)\right)}+u^{1/2}g\left(\tr(X)\right)\right)^{-(2^{j+1}-1)}}
\left(\sum_{k=0}^{\frac{n-1}{2}}X^{q^{2k}}\right)^{2^j}\right].
\end{aligned}
\end{equation*}
This yields $\tilde{F}(x)$ in Theorem \ref{recuinvconst}.

\section{Concluding remarks}\label{seccon}
In this paper, a construction of complete permutation polynomials
over finite fields of even characteristic is generalized in a
recursive manner. In addition,  several  classes of complete
permutation polynomials are derived by inverting known ones. In
fact, the idea of seeking a recursive construction of complete
permutation polynomials is motivated by the work in \cite{lc1} on
constructing bilinear permutation polynomials over finite fields of
even characteristic. It should be noted that from our construction,
complete permutation polynomials over any extension field of
$\mathbb{F}_2$ can be obtained provided the extension degree has an
odd factor. It is a natural question that how to construct complete
permutation polynomials over ramified extensions of $\mathbb{F}_2$
(i.e., the extension degrees are powers of two). To the best of the
authors' knowledge, no special attention has been paid to this
question before. We leave it as an open problem.\medskip

\textbf{Open problem.} Construct complete permutation polynomials
over the finite field $\mathbb{F}_{2^{2^e}}$, where $e$ is a
positive integer.

\section*{Acknowledgements}
The authors would like to thank Prof. Michael Zieve for helpful
comments on an earlier version of the paper. This work is partially
supported by the National Key Basic Research Project  of China under
Grant No. 2011CB302400, the National Natural Science Foundation of
China under Grant No. 61379139, and the ``Strategic Priority
Research Program" of the Chinese Academy of Sciences under Grant No.
XDA06010701.


\end{document}